\documentclass{amsart}
\usepackage{amssymb,amscd}
\usepackage{xspace}

\newtheorem{theo}{Theorem}
\newtheorem{prop}[theo]{Proposition}

\def\algc{alg\-ebr\-aic\-ally closed}
\def\char{characteristic}
\def\d{{\delta}}
\def\iff{if and only if} 
\def\P{{\bf P}}

\def\Q{{\bf Q}}

\def\cO{{\mathcal O}}
\def\cE{{\mathcal E}}

\begin{document}
\title{Fano varieties with high degree}

\author{Olivier Debarre}
\address{Institut de Recherche Math\'ematique Avanc\'ee\\
Universit\'e Louis Pasteur et CNRS\\ 7 rue Ren\'e Descartes\\ 67084
Strasbourg C\'edex\\ France}
\email{debarre@irma.u-strasbg.fr}
\urladdr{http://www-irma.u-strasbg.fr/\~{}debarre/}

\maketitle

A Fano variety is a smooth projective variety defined over an
\algc\ field.  

By Yau's proof of   Calabi's conjecture (\cite{yau1},
\cite{yau2}), complex Fano varieties are compact K\"ahler varieties
with positive definite Ricci curvature (\cite{Bou},
\cite{besse}, 11.16.ii)). In particular, compact complex varieties
with a positive K\"ahler-Einstein metric are Fano varieties. 

Not all complex Fano varieties however have a K\"ahler-Einstein
metric, because this forces the automorphism group  to be reductive 
(\cite{besse}, cor.~11.54).  This excludes for example the
projective plane blown-up at one or two points.  However, again,
this condition is not sufficient to ensure the existence of  a
K\"ahler-Einstein metric (Tian gives in \cite{tian} an example of a
Fano variety with no K\"ahler-Einstein metric and finite
automorphism  group). The existence of such a metric also implies
that the tangent bundle is
$(-K_X)$-stable.

If $X$ is   a Fano variety of dimension $n$,   we denote by $\d(X)$
the $n$th root  of the intersection number  $ (-K_X)^n$. We let
$\rho_X$ be the rank of the N\'eron-Severi group of $X$ (also called
the Picard number of $X$), and $\iota_X$ the {\em index}  of $X$,
i.e. the greatest integer by which the  canonical class is
divisible; it satisfies $\iota_X\le n+1$, with equality \iff\ $X$ is
isomorphic to $\P^n$.  Various upper bounds on $\d(X)$ are known
when the base field has
\char\ zero.
\begin{itemize}
\item For any Fano variety $X$ of dimension $n$, we have
(\cite{KMMc})
$$\d(X)\le 3(2^n-1)(n+1)^{(n+1)(2^n-1)}\ .$$
\item When $\rho_X=1$, we have (\cite{ca}, \cite{na},
\cite{KMMa}, \cite{ran})
$$\d(X)\le \max(n\iota_X,n+1)\le n(n+1)\ .$$
\item When $\rho_X=1$ and the tangent bundle of $X$ is
$(-K_X)$-semi-stable, we have (\cite{ran})
$$\d(X)\le 2n\ .$$
\end{itemize} When $X$ has a K\"ahler-Einstein metric, classical
methods of differential geometry give (see \cite{deb})  
$$\d(X)\le (2n-1)\Biggl( \frac{2^{n+1}(n!)^2}{ (2n)!}\Biggr)^{1/n}
\sim 2n$$ which is asymptotically the same as Ran's bound. Finally,
Reid proved in \cite{reid} that the tangent bundle of a Fano variety
$X$ with
$\rho_X=\iota_X=1$ is $(-K_X)$-stable.
\bigskip

{\em The purpose of this note is to construct, for each positive
integers
$k$ and $n$, Fano varieties of dimension $n$ and Picard Number $k$
for which $\d(X)$ grows essentially like $n^k$.}

\bigskip
 Batyrev remarked in
\cite{batyrev} that for the $n$-dimensional Fano variety\footnote{We
follow Grothendieck's notation: for a vector bundle
$\cE$, the projectivization $\P\cE$ is the space of {\em
hyperplanes} in the fibers of $\cE$.} 
$$X=\P
\bigl(\cO_{\P^{n-1}}\oplus
\cO_{\P^{n-1}}(n-1)\bigr)\ ,$$  we have 
$$\d(X) =\Biggl( \frac{(2n-1)^n-1}{ n-1} \Biggr)^{1/n} \sim 2n\ .$$
Consider more generally
$$X =\P(\cO_{\P^s}^{\oplus r}\oplus\cO_{\P^s}(a))\ ,
$$  where $r$, $s$ and $a$ are non-negative integers. We have
$$-K_X\sim (r+1)L+ (s+1-a)H \ ,$$ where $L$ is a divisor associated
with the line bundle
$\cO_X(1)$ and $H$ is the pull-back on $X$ of a hyperplane in
${\P^s}$. It follows that $X$ is a Fano variety when $a\le
s$~\footnote{This can   be seen  by noting that $\cO_X(L+H)$ is the
line bundle $\cO(1)$ on $X$ associated with the description of $X$
as $ \P(\cO_{\P^s}(1)^{\oplus r}\oplus\cO_{\P^s}(a+1))$. It is
therefore ample.}. 

   In the intersection ring of $X$, we have the relations
$L^{r+1}=a H\cdot L^r$ and $L^r\cdot H^s=1$. Setting
$n=\dim(X)=r+s$, we get
\begin{align*} (-K_X)^n&= \sum_{i=0}^n\binom{n}{i}(r+1)^iL^i\cdot
H^{n-i}\\ &=\sum_{i=r}^n\binom{n}{i}(r+1)^i(aH)^{i-r}\cdot L^r\cdot
\pi^*H^{n-i}\\ &=\sum_{i=r}^n\binom{n}{i}(r+1)^ia^{i-r}\\ &\ge
(r+1)^na^{n-r} 
\end{align*} Take $a=s=n-r$; the function $r\mapsto r^n(n-r)^{n-r}$
reaches its maximum near $\frac{n}{\log n}$. Taking
$r=\bigl[\frac{n}{\log n}\bigr]$, we get
\begin{align*} (-K_X)^n&\ge
\Bigl( \frac{n}{\log n}\Bigr)^n
\Bigl(n-\frac{n}{\log n}\Bigr)^{n-\frac{n}{\log n}}\\ &\ge
n^{2n-\frac{n}{\log n}}
\ \frac{1}{(\log n)^n}
\ \Bigl(1-\frac{1}{\log n}\Bigr)^n\\ &= n^{2n}e^{-n}
\ \frac{1}{(\log n)^n}
\ \Bigl(1-\frac{1}{\log n}\Bigr)^n\\ &\ge \Bigl( \frac{3n^2}{10\log
n}\Bigr)^n
\end{align*} for $\log n\ge \frac{10}{10-3e}$, i.e. $n\ge 226$. This
lower bound for $ (-K_X)^n$   actually holds for all
$n\ge 3$ by direct calculation. Furthermore, even when taking the
value of $r$ which gives the highest degree, numerical calculations
show  that
$\d(X)$ is still equivalent to   some (non-zero) multiple of 
$\frac{n^2}{\log n}$.

\begin{prop} For each $n\ge 3$, there is a Fano variety $X$ of
dimension $n$, index $1$ and Picard number $2$ such that
$$\d(X)\ge  \frac{3n^2}{10\log n}\ .$$
\end{prop}

\medskip If we analyze the construction, we see that we need to 
start  from a Fano variety with both high index and high degree. The
variety $X$ constructed in the proposition has index $1$, hence
cannot be used to iterate the process. However, if one takes instead
$a=s-r=n-2r$ and the same $r$, the index of $X$ becomes $r+1$, and,
although the degree becomes slightly smaller, it still satisfies
\begin{align*} (-K_X)^n&\ge
\Bigl( \frac{n}{\log n}\Bigr)^n
\Bigl(n-2\frac{n}{\log n}\Bigr)^{n-\frac{n}{\log n}}\\ &\ge \Bigl(
\frac{n^2}{7\log n}\Bigr)^n
\end{align*} for $\log n\ge \frac{14}{7-e}$; this lower bound  
actually holds for all
$n\ge 4$.

\begin{prop} For each positive integers $k\ge 2$ and $n\ge 4$ such
that 
$\frac{n}{\log n}\ge 2^{k-2}$, there exist  a positive   constant
$c(k)$~\footnote{One can take $c(k)=
\displaystyle\frac{1}{4^{k^2-k+2}}$.} and a Fano variety
$X$ of dimension $n$ and Picard number $k$ such that
$$ 
\d(X)\ge \frac{c(k)n^k }{ (\log n)^{k-1}} \ .
$$\end{prop}

\begin{proof} We proceed by induction on $k$, {\em assuming in
addition that the index of $X$ is} $\Bigl[\frac{n}{2^{k-2}\log
n}\Bigr]+1$. We just did it for $k=2$. Assume the construction is
done for $k\ge 2$. Let $n$ be an integer as in the proposition, and
set
$$r=\Bigl[\frac{n}{2^{k-1}\log n}\Bigr]   \hskip 15mm  s=n-r\ .
$$ Since
$\frac{n}{\log n}\ge 2^{k-1}$, the integer $r$ is positive. Also,
$r\le\frac{n}{4}$ because $n\ge 4$ and $k\ge 2$. It implies
$$\frac{s}{\log s}\ge \frac{3n}{4\log n}> 2^{k-2}\ ,$$ hence there
exists by induction a Fano variety
$Y$ of dimension $s$, index
 $\iota_Y=\Bigl[\frac{s}{2^{k-2}\log s}\Bigr]+1$ and Picard number
$k$, such that 
$$\d(Y)\ge \frac{c(k)s^k}{( \log s)^{k-1}} \ .
$$ Write
$-K_Y=\iota_Y H$, with $H$ ample on $Y$, and let 
$$X =\P(\cO_Y^{\oplus r}\oplus\cO_Y((\iota_Y-r-1)H))\ ,$$ with
projection $\pi:X\to Y$, so that $-K_X=(r+1)(L+\pi^*H)$. As above,
it implies that
$X$ is a Fano variety of dimension $n=r+s$ {\em when} $r<
\iota_Y 
$; note also that $\rho_X=k+1$ when $r>0$, and
$\iota_X=r+1$. We get again
\begin{align*} (-K_X)^n&= (r+1)^n\sum_{i=0}^n\binom{n}{i}L^i\cdot
\pi^*H^{n-i}\\
&=(r+1)^n\sum_{i=r}^n\binom{n}{i}((\iota_Y-r-1)H)^{i-r}\cdot L^r\cdot
\pi^*H^{n-i}\\
&=(r+1)^n\sum_{i=r}^n\binom{n}{i}(\iota_Y-r-1)^{i-r}H^s\\ & \ge
(r+1)^n(\binom{n}{r}+(\iota_Y-r-1)^s)H^s\\ &\ge
\Bigl(\frac{n}{2^{k-1}\log n}\Bigr)^n   (1+(\iota_Y-r-1)^s)
\Biggl(\frac{c(k)s^k}{( \log s)^{k-1}}\Biggr)^s
\frac{1}{\iota_Y^s}
\end{align*} Note that
$$
\iota_Y\ge \frac{s}{2^{k-2}\log n}\ge \frac{3n}{2^k\log n}\ .
$$ If $\frac{n}{\log n}\ge 7\cdot 2^{k-2}$, we obtain
$$\frac{r+1}{\iota_Y}\le \frac{2}{3}+\frac{2^k\log n}{3n}\le
\frac{6}{7}\ ;
$$ if $\frac{n}{\log n}< 7\cdot 2^{k-2}$, we get
$$\frac{1}{\iota_Y}\ge \frac{1}{\frac{s}{2^{k-2}\log s}+1}\ge
\frac{1}{\frac{n}{2^{k-2}\log n}+1}\ge\frac{1}{8}\ .
$$ In all cases, 
$$(1+(\iota_Y-r-1)^s)\frac{1}{\iota_Y^s}\ge\frac{1}{8^n}\ .
$$ It follows that
\begin{align*}
\d(X) &\ge \frac{n}{2^{k-1}\log n}
\Bigl(\frac{c(k)(3n)^k}{4^k( \log n)^{k-1}}\Bigr)^{1-\frac{r}{n}}
\ \frac{1}{8}\\ &\ge n^{1+k-k\frac{1}{2^{k-1}\log n}}
\frac{1}{2^{k+2}\log n}
\ \frac{c(k)3^k}{4^k( \log n)^{k-1}}\\ &\ge n^{k+1}\ \frac{1}{(\log
n)^k}\ 
\frac{c(k)3^k}{e2^{k+2}4^k }\ ,
\end{align*} which proves the proposition.
\end{proof}

\end{document}